\documentclass[a4paper,11pt]{article}

	\usepackage{amsfonts}\usepackage{amsthm}\usepackage{amsmath} \usepackage{a4}
	\newcommand{\R}{\mathbb R}

	\newcommand{\cP}{{\mathcal P}} \newcommand{\cT}{{\mathcal T}} 
	\newcommand{\cV}{{\mathcal V}} \newcommand{\cA}{{\mathcal A}}
	
	\newcommand{\cQ}{{\mathcal Q}} \newcommand{\cS}{{\mathcal S}}
    \newcommand{\cR}{{\mathcal R}}
 	
 	\newcommand{\bd}{\operatorname{\partial}}
    \newcommand{\und}{\,\,\,\hbox{ and } \,\,\,} \newcommand{\fuer}{\,\,\,\hbox{ for } \,} 

	\newcommand{\eqnref}[1]{(\ref{#1})} 

	\newtheorem{theorem}{Theorem}  
	\newtheorem{lemma}{Lemma}  
	\newtheorem*{proposition}{Proposition}
	\newtheorem{corollar}[theorem]{Corollary}

\title{Elementary moves on triangulations}

\author{Monika Ludwig and Matthias Reitzner\footnote{Research of both authors was supported in part by the European Network PHD, MCRN-511953.}}

\date{}

\begin{document}

\maketitle
 \begin{abstract}
It is proved that a triangulation of a polyhedron can always be
transformed into any other triangulation of the polyhedron 
using only elementary moves. One consequence is that an additive
function (valuation) defined only on simplices may always be extended to an
additive function on
all polyhedra.
\bigskip
 
 
\end{abstract}

An $n$-polyhedron $P$ in $\R^{N}$, $1\le n \le N$, is a finite union of
$n$-dimensional polytopes, where a polytope is the compact convex hull of finitely many points in $\R^N$.  A finite set of $n$-simplices $\alpha P$
is a {\em triangulation} of $P$ if no pair of simplices intersects in a set of dimension $n$ and their union equals $P$.
We shall investigate transformations of triangulations by elementary
moves. Here an {\em elementary
move\,} applied to $\alpha P$ is one of the two following operations:
a simplex $T\in\alpha P$ is dissected into two $n$-simplices $T_1, T_2$
by a hyperplane containing an
$(n-2)$-dimensional face of $T$; or the
reverse, that is, two simplices $T_1, T_2\in\alpha P$ are replaced by
$T=T_1\cup T_2$ if $T$ is again a simplex.
We say that triangulations $\alpha P$ and $\beta P$ are equivalent by
elementary moves,
and write $\alpha P \sim \beta P$, if there  are finitely many elementary
moves that transform $\alpha P$ into $\beta P$. The main object of this note
is to show the following result.

\begin{theorem}\label{move}
If  $\,\alpha P$ and  $\beta P$ are triangulations of the $n$-polyhedron
$P$, then $\alpha P\sim \beta P$.
\end{theorem}

A triangulation  $\alpha P$ with the additional property that any pair of
simplices intersects in a common face gives rise to a simplicial complex
$\hat\alpha P$. 
It is a classical result of algebraic topology due to
Alexander \cite{Alexander30} and Newman \cite{Newman29, Newman31} (see also  \cite{Lickorish99}) that
a simplicial complex $\hat \alpha P$ can always be
transformed into any other simplicial complex $\hat \beta P$ with the same underlying polyhedron by using only finitely many  stellar moves.
Here a {\em stellar move} is a suitable sequence of elementary
moves followed by a simplicial isomorphism, where 
a {\em simplicial isomorphism} between two complexes is a bijection between
their vertices that induces a bijection between their $k$-dimensional simplices for $1\le k\le n$. For  precise definitions, see  \cite{Lickorish99} and for related results, see \cite{Ewald:Shephard}, \cite{Pachner81}, \cite{Pachner90}.
The new feature of
Theorem~\ref{move} is that simplicial isomorphisms are not allowed.
So our theorem belongs to metric geometry whereas the Alexander-Newman
theorem is a topological result.

As an application of Theorem \ref{move} we obtain the following results on
valuations.  Here a function $\mu: \cS\to \R$ defined on a class $\cS$ of sets is called {\em a valuation} or {\em additive} if $\mu(\emptyset)=0$, where $\emptyset$ is the empty set, and if
\begin{equation*}
\mu(S) +  \mu(T) = \mu(S\cup T) + \mu(S\cap T),
\end{equation*}
for  all $S,T  \in \cS$  such that $S\cup T, S\cap T \in {\cS}$ as well. 
\goodbreak

Let $\cT^n$ be the set of simplices of dimension at most $n$ in $\R^N$.
Note the following connection between the definition of valuations on simplices and elementary moves.
If $S,T\in\cT^n$, then $S\cup T\in\cT^n$ implies that $S$ and $T$ can be obtained from $S\cup T$ by elementary moves and the most basic case is when $S\cup T$ is dissected by an elementary move into $S$ and $T$. Let $\cQ^n$ be the set of polyhedra of dimension at most $n$ in $\R^N$

\begin{theorem}\label{extend}
Every valuation on $\cT^n$ has a unique extension to a valuation on $\cQ^n$.
\end{theorem}
 
\noindent
As a corollary we obtain the analogous theorem for polytopes. Let $\cP^n$ be the set of polytopes of dimension at most $n$ in $\R^N$.

\begin{corollar}\label{extendcp}
Every valuation on $\cT^n$ has a unique extension to a valuation on $\cP^n$.
\end{corollar}

\noindent
A version of Corollary \ref{extendcp} is stated and used in \cite{Kuperberg03}. Note that in the proof of this result in \cite{Kuperberg03} the Alexander-Newman theorem has to be replaced by  Theorem \ref{move} of the present paper. 
  
Valuations on polyhedra are classical and they played a critical role
in Dehn's solution of  Hilbert's Third Problem. Results regarding the
classification and characterization of invariant valuations are central to
convex and integral geometry; see \cite{Hadwiger:V}, \cite{Klain:Rota},
\cite{McMullen93}, \cite{McMullen:Schneider}. In recent years, many new
results on valuations have been obtained, see, for example,
\cite{Alesker99}--\nocite{Alesker01}\nocite{Alesker03}\nocite{Alesker04}\cite{Alesker04},
\cite{Klain95}--\nocite{Klain96}\nocite{Klain97}\cite{Klain00},
\cite{Ludwig:origin}--\nocite{Ludwig:matrix}\cite{Ludwig:Reitzner},
\cite{Schneider96}.  
Also extension questions for general valuations are classical. Volland
\cite{Volland} and Perles and Sallee \cite{Perles:Sallee} proved that every
valuation on polytopes has a unique extension to a valuation on
polyhedra. Their results were generalized by Groemer \cite{Groemer78}.
Volume is the most basic example of a valuation on polyhedra. The question
of how to extend volume from simplices  to polyhedra is closely
connected with the
quest for an elementary definition of volume for polytopes.
In his Third Problem, Hilbert  \cite{Hilbert00} asked whether volume on polyhedra can be defined
by using only scissors congruences, that is, if two polyhedra
$P$ and $Q$ of the same volume can each be cut into a finite number
of pieces $P_1,\ldots, P_m$ and $Q_1,\ldots,Q_m$ with $P_i$ congruent to $Q_i$
by a rigid  motion for each $i$. This question was answered in the negative by Dehn \cite{Dehn01} (see also \cite{Boltianskii, Sah}).  However,  there are simple definitions for volume of simplices. So the volume of $n$-simplices can be defined as  height times the
$(n-1)$-dimensional volume of its base divided by $n$. A simple geometric argument (see \cite{Suess}) shows that this does not depend on the choice of the base and volume defined in this way is a valuation on simplices.  
Schatunovsky
\cite{Schatunovsky} (for $n=3$) and S\"uss \cite{Suess} (for general
dimensions) proved that there
is a unique extension of this valuation from simplices to polyhedra. In that way they
obtained an elementary definition of volume for polyhedra. Theorem
\ref{extend} is the extension of their result to general valuations.

Theorem \ref{move} has already been used by a number of authors. In
fact, results more general than Theorem \ref{move} have already been used.
For example, Lemma 2.2 in
Sah's book on
Hilbert's Third Problem \cite{Sah} states that any
two triangulations have a common refinement by elementary moves using only
dissections. However, this is well-known to be an important open
problem in algebraic
topology, see \cite{Lickorish99}. Later, Sah \cite{Sah83} replaced his lemma
by Theorem \ref{move} of the present paper. While Sah observes that Theorem \ref{move}
suffices for the
constructions in his book, he neglects to provide a proof of  Theorem
\ref{move}. For related
results, see also \cite{Dupont}.
 
The authors thank Professors G.~Kuperberg, W.B.R.~Lickorish, and
P.~McMullen for help in clarifying which statements and proofs are
available in the literature and  Professors P.M.~Gruber and  E.~Lutwak  and the anonymous referees for their helpful comments.

\section{Proof of Theorem \ref{move}}

A triangulation $\alpha P$ of an $n$-polytope $P$ is called a {\em starring} at $a\in P$, if every $n$-simplex in $\alpha P$ has a vertex at $a$. Note that every $n$-polytope $P$ has a starring at every $a\in P$. This can be seen by using induction on $n$. It is trivial for $n=1$. Suppose that  every  $(n-1)$-polytope has a starring. Let $P$ be an $n$-polytope and $a\in P$.  For every $(n-1)$-dimensional face $F_j$ of $P$ with $a\not\in F_j$, we choose a starring $\alpha_j F_j$. Then the convex hulls of $a$ and the $(n-1)$-simplices in $\alpha_j F_j$ are a starring of $P$ at $a$.

A triangulation $\gamma P$ of an $n$-polyhedron $P$ is a {\em refinement} of $\alpha P$ if every simplex of $\gamma P$ is contained in a simplex of $\alpha P$. Note that any two triangulations $\alpha P$ and $\beta P$ have a common refinement. To see this, let $P_1, \ldots, P_l$  be the  polytopes $S\cap T$, $S\in \alpha P$, $T\in\beta P$, that are $n$-dimensional.  For $j=1,\ldots, l$, we choose a triangulation $\gamma_j P_j$  of $P_j$, for example, 
we can take a starring of $P_j$ at any point $a_j\in P_j$. Then $\gamma P=\gamma_1 P_1\cup\dots\cup \gamma_l P_l$ is a triangulation of $P$. Since every simplex in $\gamma P$ is contained in a suitable $P_j$, $\gamma P$ is a common refinement of $\alpha P$ and $\beta P$. 

Let $P$ be an $n$-polyhedron. To show that $\alpha P\sim \beta P$ for any two triangulations $\alpha P$ and $\beta P$ of $P$, we prove that any refinement of a given triangulation can be 
obtained from the original triangulation by finitely many elementary moves. 
To show this, it is enough to prove the following proposition.

\begin{proposition}
If $\alpha T$ is a triangulation of the $n$-simplex $T$, then $\alpha T\sim T$.
\end{proposition}

\noindent
Here we write $T$ (instead of $\{T\}$) for the trivial triangulation of $T$.

The rest of this section is devoted to the proof of this proposition. We follow the classical approach of Alexander and Newman as presented by Lickorish \cite{Lickorish99}. The main new feature of the proof is Lemma \ref{our}.

We use induction on the dimension $n$. The case $n=1 $ is trivial. Assume that the proposition is true for dimensions less or equal $n$, that is, for every triangulation $\alpha S$ of a simplex $S$
\begin{equation}\label{inddim}
\alpha S \sim S \fuer \dim S<n,
\end{equation}
where $\dim$ stands for dimension.

Note that, if $S$ is a simplex,  then every elementary subdivision of $S$ induces an elementary subdivision of the simplices in the boundary of $S$.
The following observation is used several times.
Let $P= [Q,v]$ be a pyramid with apex $v$ and base $Q$. Let $\alpha Q=\{S_1,\ldots, S_k\}$ and $\beta Q=\{ T_1,\ldots, T_l\}$ be triangulations of the $(n-1)$-polytope $Q$.
Let $[\alpha Q,v] =\{[S_1,v],\ldots,[S_k,v]\}$ and $[\beta Q,v]=\{[T_1,v],\ldots,[T_l,v]\}$. Then  $[\alpha Q, v]$ and $[\beta Q, v]$ are  triangulations of $P$ and
\begin{equation}\label{pyramid}
\alpha Q\sim \beta Q \,\,\,\Rightarrow\,\,\,  [\alpha Q, v] \sim [\beta Q, v].
\end{equation}

\begin{lemma}\label{Newman}
If  $\,\alpha T$ is a starring of an $n$-simplex $T$, then $\alpha T\sim T$.
\end{lemma}

\begin{proof}
Assume $T\subset \R^n$.
First, let  $T$ be starred at $a\in \bd T$. Let $S_0, \ldots, S_n$ be the facets of $T$ and let $a\in S_0$. Then
$$\alpha T=[\alpha_1 S_1, a]\cup \dots \cup [\alpha_n S_n,a],$$
where $\alpha_i S_i$ is a triangulation of $S_i$.
The induction hypothesis \eqref{inddim} and \eqnref{pyramid} imply that 
$[\alpha_i S_i, a]\sim [S_i,a].$
Thus, 
$$\alpha T \sim [S_1, a]\cup \dots \cup [S_n,a].$$
Let $v$ be the vertex of $T$ opposite to $S_0$. Then there is a starring $\alpha_0 S_0$  of $S_0$ at $a$ such that 
$$[S_1, a]\cup \dots\cup [S_n,a] =[\alpha_0 S_0, v].$$
The induction hypothesis \eqref{inddim} and \eqnref{pyramid} imply that $[\alpha_0 S_0, v] \sim [S_0,v]$. Consequently,
\begin{equation}\label{fc}
\alpha T\sim T \fuer a\in \bd T.
\end{equation}

Now, let $a$ be a point in the interior of $T$. Then
$$\alpha T = [\alpha_0 S_0, a]\cup \dots\cup [\alpha_n S_n,a],$$
where $\alpha_i S_i$ is a triangulation of the facet $S_i$ of $T$. By the induction hypothesis \eqref{inddim} and \eqnref{pyramid} 
$$\alpha T\sim [S_0, a]\cup \dots\cup [S_n,a].$$
We write $a T$ for the starring at $a$ of the $n$-simplex $T$ if every simplex in $a T$ has a facet of $T$ as its base.
We dissect $T$ by a hyperplane $H$ through $a$ and an $(n-2)$-dimensional face of $T$ into two simplices $T^+, T^-$. This is an elementary move. Thus by \eqnref{fc}
$$T\sim T^+\cup T^-\sim a T^+ \cup a T^-.$$
If a facet $S$ is subdivided by $H$ into $S^+$ and $S^-$, then $[S^+,a]\cup [S^-,a]\sim [S,a]$.
Thus
$$a T^+ \cup a T^-\sim a T.$$
This completes the proof of the lemma.\end{proof}

Note that, for $a \in P$ fixed, \eqref{inddim} and \eqref{pyramid} imply that any two starrings at $a$ are equivalent by elementary moves. We write $a P$ for a starring of a polytope $P$ at a point $a\in P$.

\begin{lemma}\label{our}
For every $n$-polytope  $P$ and $a,b\in P$, $a P\sim b P$.
\end{lemma}

\begin{proof}
Assume $P\subset \R^n$.
We use induction on the number $m$ of vertices of $P$. If  $m=n+1$, then $P$ is an $n$-simplex and the statement is true by Lemma \ref{Newman}. Suppose the lemma is true for polytopes with at most $m$ vertices. Write $[A_1, \ldots, A_l]$ for the convex hull of sets $A_1, \ldots, A_l$.

Let $P$ be a polytope with vertices $v_1, \ldots, v_{m}, v$. Let $a \in P^- = [v_1, \dots ,v_{m}]$. 
We say that a facet $F$ of $P^-$ is visible from $v$ if for every $x\in F$, $[v,x]\cap P^-=\{x\}$. By starring $P^-$ at $a$, the facets of $P^-$ that are visible from $v$ are subdivided into $(n-1)$-simplices $S_i$, $i=1,\ldots,l$. Let  $\cV_a$ be the set of the $n$-simplices $[S_i,v]$, $i=1,\ldots,l$. Let $a P^-$ contain the simplices $[S_i,a]$, $i=1, \dots, l$, that are $n$-dimensional.

The main step is to show that
\begin{equation}\label{vk+}
a P \sim a P^- \cup \cV_a.
\end{equation}
If \eqnref{vk+} holds and if $b \in P^-$, then the induction hypothesis on the number of vertices implies that $a P^-\sim b P^-$. By \eqnref{inddim} and \eqnref{pyramid} we obtain $\cV_a\sim \cV_b$. Thus
$$a P \sim a P^- \cup \cV_a \sim b P^- \cup \cV_b \sim b P.$$
If there is no vertex $v$ such that $a,b\in [v_1,\ldots, v_m]$, then there are vertices $v,w$ of  $P$ such that  $a \in [v_1, \dots, v_{m-1}, w],\ b \in [v_1, \dots, v_{m-1}, v]$ and  $P = [v_1, \dots, v_{m-1}, v, w]$. Choose $c \in [v_1, \dots, v_{m-1},w] \cap [v_1, \dots, v_{m-1},v]$ in the interior of $P$. Then by \eqref{vk+}
$$a P \sim c P \und c P \sim b P.$$
Thus it is enough to show \eqnref{vk+} to prove the lemma.

Let $H_j$, $j=1, \dots, k$, be the affine hulls of $S_i$, $i=1,\ldots,l$.
Denote the intersection point of the hyperplane $H_j$ with the segment  $[a, v]$ by  $x_j$ (these points are not necessarily distinct). Without loss of generality assume that the hyperplanes $H_j$ are numbered such that $[a,v]$ is dissected into  $[a, x_1], [x_1, x_2], \ldots,[x_{k-1}, x_{k}],[x_{k}, x_{k +1} =v]$.
Let  $\cA_1$ be the set of those  simplices $[S_i,a]$, $i=1,\ldots, l$, that are $n$-dimensional. Set $\cV_1=\cV_a$.

First, assume that $a\ne x_1$ and $x_k\ne v$.
Starting with $j=1$, we construct $\cA_{j+1}$ and $\cV_{j+1}$
from $\cA_j$ and $\cV_j$ in the following way. The simplices in $\cA_{j}$ have $a$ as a vertex. Let $\cA_j (H_j)\subset \cA_j$ be the subset of those simplices that have their $(n-1)$-dimensional base in  $H_j$.
The simplices in $\cV_{j}$ have $v$ as a vertex.
Let $\cV_j (H_j)\subset \cV_j$ be the subset of  those simplices that have their $(n-1)$-dimensional base in $H_j$. Note that the set of these $(n-1)$-dimensional bases coincides for $\cA_j(H_j)$ and $\cV_j(H_j)$.
These bases form a facet $F$, say,  of $[P^-,x_j]$. By the induction hypothesis, the triangulation of $F$ by these $(n-1)$-dimensional bases is equivalent by elementary moves to a starring of $F$ at $x_j$. The moves applied in $H_j$ induce moves on $\cA_j (H_j)$ and $\cV_j (H_j)$. Denote the sets obtained by these moves by $\cA'_j(H_j)$ and $\cV'_j (H_j)$. Note that the set of $(n-1)$-dimensional bases coincides for $\cA_j' (H_j)$ and $\cV_j' (H_j)$. We replace $\cA_j (H_j)$ by $\cA_j' (H_j)$ and $\cV_j (H_j)$ by  $\cV_j' (H_j)$ and set 
$ \cA'_j=(\cA_j\backslash \cA_j (H_j))\cup \cA_{j}'(H_j)$ and $\cV_j' =(\cV_j\backslash \cV_j (H_j))\cup \cV_{j}'(H_j)$.

In the next step, we do the following. If $x_j=x_{j+1}$, we set $\cA_{j+1}=  \cA_{j}'$ and
$\cV_{j+1}=  \cV_{j}'$ and  have
$$\cA_j \cup \cV_j \sim \cA_{j+1} \cup \cV_{j+1} . $$
If $x_j \neq x_{j+1}$,  denote by $\cA_j'(x_j)\subset  \cA_j'$ the subset of  those simplices  which contain $x_j$ as a vertex and by $\cV_j'(x_j)\subset \cV_j'$ the subset of  those simplices which contain $x_j$ as a vertex. Each simplex in $\cV_j'(x_j)$ is the convex hull of $x_{j}$, $v$ and an $(n-2)$-dimensional face $B$, say. The convex hull $[B, x_j, v]$ is in $\cV_j'(x_j)$ and  $[B,a, x_j]$ in $\cA_j'(x_j)$.
We subdivide the simplex $[B, x_j, v]$ into the simplices $[B, x_j, x_{j+1}]$ and $[B, x_{j+1}, v]$.
Then we take the union of $[B, a, x_j]$  and $[B, x_j, x_{j+1}]$ and obtain the simplex $[B, a, x_{j+1}]$.
These operations are elementary moves. We do this with every simplex with vertex $x_j$. We obtain the new set  $\cA_{j}'(x_{j+1})$ which contains the  simplices $[B, a, x_{j+1}]$ and the new set $\cV_{j}'(x_{j+1})$  which contains the simplices $[B, x_{j+1}, v]$.  Let $\cA_{j+1}= (\cA_j'\backslash \cA_j' (x_j))\cup \cA_{j}'(x_{j+1})$ for $j\le k$ and
$\cV_{j+1}= (\cV_j'\backslash \cV_j' (x_j))\cup \cV_{j}'(x_{j+1})$ for $j<k$ and $\cV_{k+1}=\emptyset$. We have 
$$\cA_j \cup \cV_j \sim \cA_{j+1} \cup \cV_{j+1}. $$

Next, we consider the case $a=x_1=\dots=x_i$. Note that each facet of $P^-$ that contains $a$ is already starred at $a$.
Let $\cV_1(a) \subset \cV_1$ be the subset of  those simplices which contain $a$ as a vertex. Each simplex in $\cV_1(a)$ is the convex hull of $a$, $v$ and an $(n-2)$-dimensional face $B$. We subdivide the simplex $[B, a, v]$ into the simplices $[B, a, x_{i+1}]$ and $[B, x_{i+1}, v]$.
These operations are elementary moves. We do this with every simplex with vertex $a$. We obtain the new set  $\cA_{i+1}$ which contains the  simplices $[B, a, x_{i+1}]$ and the new set $\cV_{i+1}$  which contains the simplices $[B, x_{i+1}, v]$ and the simplices in 
$\cV_1\backslash \cV_1 (a)$. 
We have  $ \cV_1 \sim \cA_{i+1} \cup \cV_{i+1} $. Starting with $j=i+1$, we construct $\cA_{j+1}$ and $\cV_{j+1}$ by the algorithm described above.
If $x_{i'}=\dots=x_k=v$, then the above algorithm implies that $\cV_{i'}$ contains only lower dimensional simplices. Therefore we set $\cV_{i'}=\emptyset, \dots, \cV_{k+1}=\emptyset$ and $\cA_{i'+1}=\cA_{i'}, \dots , \cA_{k+1}=\cA_{i'}$.

Thus in all cases 
$$\cA_1 \cup \cV_1 \sim \dots  \sim \cA_{k+1}\cup\cV_{k+1}= \cA_{k+1}.$$
Since $a P \backslash \cA_{k+1}= a P^-\backslash \cA_1$, this proves \eqnref{vk+}.
\end{proof}

Let $T$ be an $n$-simplex and let $\alpha T$ be a triangulation of $T$. Assume $T\subset \R^n$. Let $H_1, \dots, H_l$ be the affine hulls of the facets of the simplices in $\alpha T$.
Let $\zeta_j T$ be the dissection into polytopes (in general not simplices) of $T$ by the hyperplanes $H_1, \ldots, H_j$  and let $\zeta_0 T = T$. We dissect each polytope $P$ of $\zeta_j T$ into simplices by starring $P$ at an arbitrary interior point.
Let  $\beta_j T$ be a triangulation of $T$ obtained from  $\zeta_j T$ in this way.

We use induction on the number $k$ of hyperplanes. By Lemma \ref{Newman}, we have $T\sim \beta_0 T$. For $k\ge1$ assume that
\begin{equation}\label{ind}
T\sim \beta_j T \fuer j < k.
\end{equation}
Note that we obtain $\zeta_{k} T$ from $\zeta_{k-1} T$ by cutting by $H_k$. Every cell of $\zeta_{k-1}T$ is either unchanged or it is cut into two pieces. So let  $P \in \zeta_{k-1} T$ be cut into pieces $P_1,P_2$. We show that for $a\in P$ and $a_i\in P_i$, $i=1,2$,
\begin{equation}\label{cut}
a P \sim a_1 P_1 \cup a_2 P_2.
\end{equation}
\medskip
By Lemma \ref{our}, $a P \sim b P$ for every $b\in H_k\cap P$. Again by Lemma \ref{our}, $a_1 P_1 \sim b P_1$ and $a_2 P_2 \sim b P$.
Since any two starrings of $P$ at $a\in P$ are equivalent by elementary moves
this implies \eqref{cut}.
By our definition of $\beta_{k-1} T$ and $\beta_k T$, \eqnref{cut} implies that
$$\beta_{k-1} T\sim \beta_k T.$$
Thus, by induction,
\begin{equation}\label{grob}
T \sim \beta_0 T \sim \cdots  \sim \beta_lT.
\end{equation}
Let $\alpha T$ consist of the simplices $S_1, \ldots, S_m$.
Let $\zeta_j S_i$ be the dissection into polytopes of $S_i$ by the hyperplanes $H_1, \ldots, H_j$.
Let $\beta_j S_i$ be a triangulation that is obtained by starring each cell of $\beta_j S_i$ at a point in that cell. By Lemma \ref{our} these triangulations for different starring points are equivalent. Thus we obtain as before that
$$S_i \sim \beta_0 S_i \sim \cdots  \sim \beta_l S_i.$$
Consequently,
$$\alpha T=(S_1\cup\dots\cup S_m) \sim \cdots \sim (\beta_l S_1\cup \dots\cup \beta_l S_m) = \beta_l T.$$
Combined with \eqref{grob} this completes
the proof of the proposition.

\section{Proof of Theorem \ref{extend}}

We use induction on $n$. The case $n=0$ is straightforward. Suppose that every valuation on $\cT^{n-1}$ has a unique extension to a valuation on $\cQ^{n-1}$. 
Set $\cR^n=\{ T \cup Q\, |\, T \in \cT^n,\, Q \in  \cQ^{n-1} \}$. By the induction hypothesis, it is straightforward that $\mu$ has a unique extension from $\cT^n$ to $\cR^n$. 
We prove that $\mu$ has a unique extension from $\cR^n$ to a valuation on $\cQ^n$.

Let $Q\in\cQ^n$ and let $\delta Q=\{R_1,\ldots,R_k\}$, $R_i\in \cR^n$, be a dissection  of $Q$; that is,
the intersection of a pair of elements of $\delta Q$ has dimension less than $n$  and their union equals $Q$.
For an ordered $j$-tuple $I=\{i_1,\ldots,i_j\}$, $1\le i_1<\cdots <i_j\le k $, $1\le j \le k$, set $R_I= R_{i_1}\cap \cdots \cap R_{i_j}$ and $|I|=j$. The set $\cQ^n$ is a lattice. Thus,
if $\mu$ can be extended to a valuation on $\cQ^n$, then by the
inclusion-exclusion principle
\begin{equation}\label{inex}
\mu(Q)= \sum_{I\subset\{1,\dots, k\}} (-1)^{|I|-1} \mu(R_I).
\end{equation}
We show that, if $Q\in\cQ^n$ and if $\delta Q=\{R_1,\ldots,R_k\}$, $R_i\in \cR^n$, is a dissection of $Q$, then \eqnref{inex} can be used as a definition of $\mu(Q)$.
Note that $R_I\in \cR^{n}$. Therefore the induction hypothesis implies that the right side of \eqnref{inex} is well defined. 

We show that $\mu(Q)$ as defined by \eqnref{inex} does not depend on the choice of $\delta Q$.
First, let $Q$ be an $n$-polyhedron and let $\alpha Q=\{T_1,\ldots,T_k\}$ be a triangulation of $Q$. By Theorem \ref{move}, all triangulations of an $n$-polyhedron are equivalent by elementary moves. Thus it is sufficient to show that applying an elementary move to $\alpha Q$ does not change $\mu(Q)$.
Since an elementary move is either the dissection of a simplex or the reverse, it suffices to show the following. If  $S_i=T_i$, $i=1,\ldots, k-1,$ and if $T_k$ is subdivided into $S_k$ and $S_{k+1}$, then
\begin{equation}\label{mt}
\sum_{I\subset\{1,\dots, k\}} (-1)^{|I|-1}\mu(T_I)=\sum_{J\subset\{1,\dots, k+1\}} (-1)^{|J|-1}\mu(S_J).
\end{equation}
Set $I'=I\backslash\{k\}$. Since $\mu$ is a valuation on $\cR^n$, we have for an ordered $j$-tuple $I$ with $k\in I$
\begin{eqnarray*}
\mu(T_I)&=&\mu(T_{I'}\cap T_k)\,\,=\,\,\mu (S_{I'}\cap (S_k\cup S_{k+1}))\,\,=\,\,\mu((S_{I'}\cap S_k)\cup (S_{I'}\cap S_{k+1})) \nonumber\\
&=&\mu(S_{I'}\cap S_k) +\mu(S_{I'}\cap S_{k+1})-\mu(S_{I'}\cap S_k \cap S_{k+1})\\
&=&\mu(S_{I}) +\mu(S_{I'}\cap S_{k+1})-\mu(S_{I}\cap S_{k+1}).\nonumber
\end{eqnarray*}
It follows that
\begin{eqnarray*}
\sum_{J} (-1)^{|J|-1}\mu(S_J)
&=&\sum_{k, k+1\not\in J} (-1)^{|J|-1}\mu(S_J)+\sum_{k\in J, k+1\not\in J} (-1)^{|J|-1}\mu(S_J)\\
&&+\sum_{k\not\in J, k+1 \in J} (-1)^{|J|-1}\mu(S_J)+\sum_{k,k+1 \in J} (-1)^{|J|-1}\mu(S_J)\\
&=&\sum_{k\not\in I} (-1)^{|I|-1}\mu(S_I) + \sum_{k\in I} (-1)^{|I|-1} \mu(S_{I})\\
&& +\sum_{k\in I} (-1)^{|I|-1}\mu(S_{I'}\cap S_{k+1})-\sum_{k\in I} (-1)^{|I|-1}\mu(S_{I}\cap S_{k+1})\\
&=&\sum_{k\not\in I} (-1)^{|I|-1}\mu(T_I)+\sum_{k\in I} (-1)^{|I|-1}\mu(T_I)
\,\,=\,\,\sum_{I} (-1)^{|I|-1}\mu(T_I).
\end{eqnarray*}
Thus \eqnref{mt} holds and there is a unique extension of $\mu$ to the set of $n$-polyhedra. 

Next, let $P\in\cQ^n$ be decomposed into $Q$ and $R$, where $Q$ is a uniquely defined $n$-polyhedron (or the empty set) and $R\in\cQ^{n-1}$. Let $\delta P= \alpha Q\cup R$, where $\alpha Q=\{R_1,\ldots,R_{k-1}\}$ is a triangulation of $Q$, and set $R_k =R$. This is a dissection of $P$. By \eqnref{inex},
\begin{eqnarray*}
\mu(P)&=& \sum_I (-1)^{|I|-1} \mu(R_I)
=\sum_{k\not\in I} (-1)^{|I|-1} \mu(R_I)+\sum_{k\in I}(-1)^{|I|-1} \mu(R_I)
\\ &=&\mu(Q) + \mu(R)-\mu(Q\cap R).
\end{eqnarray*}
Since $Q$ is an $n$-polyhedron and since $R,Q\cap R\in \cQ^{n-1}$, the terms $\mu(Q), \mu(R), \mu(Q\cap R)$ are well defined. Hence $\mu(P)$ does not depend on $\delta P$ either.

Finally, we show that $\mu$ as defined by \eqnref{inex} is a valuation on $\cQ^n$. Let $P,Q\in\cQ^n$.
We choose a dissection $\{R_1,\ldots,R_m\}$ of $P\cup Q$ such that, for every $R_i$,  
we have $R_i\subset P$ or $R_i\subset Q$. Then for every ordered $j$-tuple $I$,
$$\mu(P\cap R_I) +  \mu(Q\cap R_I) = \mu((P\cup Q)\cap R_I) + \mu((P\cap Q)\cap R_I).$$
By \eqnref{inex}, it follows that $\mu$ is a valuation.

\section{Proof of Corollary \ref{extendcp}}

The proof presented here relies essentially on a theorem of Tverberg \cite{Tverberg74}.  Note that it is also possible to prove the corollary by the extension theorems of  Volland \cite{Volland} or Perles and Sallee \cite{Perles:Sallee} or Groemer  \cite{Groemer78}.

The main step is to prove that there is at most one extension. 
Let $\mu_1$ and  $\mu_2$ be valuations on $\cP^n$ such that $ \mu_1 = \mu_2$ on $\cT^n$. Using induction on $n$ we show that $\mu_1 = \mu_2$ also on $\cP^n$.
The cases $n=0,1$ are trivial. Suppose that $ \mu_1 = \mu_2$ on $\cP^{n-1}$.

A {\em binary space partition} is formed by partitioning $\R^N$ by a hyperplane $H$ into two closed halfspaces $H^+$, $H^-$, and then recursively partitioning each of the two resulting halfspaces; the result is a hierarchical decomposition of space into closed convex cells (cf.~\cite{deBerg}).
Let $P$ be an $n$-polytope. Since $\mu_1$ and $\mu_2$ are valuations, after the first step of a binary space partition we have for $i=1,2,$
$$ \mu_i(P) = \mu_i(P\cap H^+)+ \mu_i(P\cap H^-) -  \mu_i(P\cap H).$$
Note that $P\cap H \in \cP^{n-1}$ for $P \not\subset H$ and thus 
$ \mu_1(P\cap H)=\mu_2(P\cap H)$. Hence to prove $ \mu_1(P) = \mu_2(P)$ it suffices to prove $\mu_1 (Q)=\mu_2(Q)$ for $Q$ equal to $P \cap H^+$ and $P \cap H^-$. 
In the next step of the binary space partition the convex polytopes $P \cap H^+$ and $P \cap H^-$ are dissected by suitable hyperplanes and the same argument applies, etc.

Tverberg \cite{Tverberg74} showed that, given an $n$-polytope $P$, one can find a binary space partition that decomposes $P$ in finitely many steps into $n$-simplices.
Since $ \mu_1 = \mu_2$ on $\cT^n$, this implies 
$\mu_1(P)=\mu_2(P)$ and the uniqueness is established.

By Theorem \ref{extend} there is at least one extension of $\mu$ to a valuation on $\cP^n$, the one which can further be extended to $\cQ^n$. Combined with the uniqueness result this proves the corollary.

\goodbreak

\bigskip\bigskip
\parindent=0pt
\begin{samepage}
Institut f\"ur Diskrete Mathematik und Geometrie\\
Technische Universit\"at Wien\\
Wiedner Hauptstra\ss e 8-10/104\\
1040 Wien, Austria\\
E-mail: monika.ludwig@tuwien.ac.at; matthias.reitzner@tuwien.ac.at
\end{samepage}

\end{document}